# A high-precision method for general nonlinear initial-boundary value problems

Jizeng Wang[1], Lei Zhang, You-He Zhou

*Key Laboratory of Mechanics on Disaster and Environment in Western China, Ministry of Education, College of Civil Engineering and Mechanics, Lanzhou University, Gansu 730000,China*

**Abstract.** A high-precision, and space–time fully decoupled, wavelet formulation numerical method is developed for a class of nonlinear initial-boundary value problems. This method is established based on a proposed Coiflet-based approximation scheme with an adjustable high order for a square-integrable function over a bounded interval, which allows expansion coefficients to be explicitly expressed by function values at a series of single points. In applying the solution method, the nonlinear initial-boundary value problems are first spatially discretized into a nonlinear initial value problem by combining the proposed wavelet approximation scheme and the conventional Galerkin method. A novel high order step-by-step time integrating approach is then developed for the resulting nonlinear initial value problem using the same function approximation scheme based on wavelet theory. The solution method is shown to have $N$th-order accuracy, as long as the Coiflet with [0, 3$N$-1] compact support is adopted, where $N$ can be any positive even number. In addition, the stability property of the method is analyzed, and the stable domain is determined. Numerical examples are considered to justify both the accuracy and efficiency of the method. Results show that the proposed solution method has better accuracy and efficiency than most other methods.

**1. Introduction**

Numerous physical phenomena, such as sound, heat, electrostatics, electrodynamics, fluid flow, elasticity or quantum mechanics, can be quantitatively described by nonlinear partial differential equations (PDEs), together with a set of additional constraints of initial, and boundary conditions, which form initial-boundary value problems (IBVPs). Given that they have an important role in science and engineering, solution methods of IBVPs have been extensively studied.

Since its introduction by Grossmann and Morlet (1984), wavelet theory has been broadly applied in many areas of science and engineering [1]. One of its successful uses is the numerical solution of PDEs because of the many attractive mathematical properties that wavelets possess [2]. Most of these studies are devoted to develop space discretization methods for nonlinear boundary value problems (BVPs) [3]. For example, a modified wavelet Galerkin method was recently proposed to solve nonlinear BVPs, including the large deflection bending of circular plates [4], rectangular plates [5], and Bratu equations [3]. A comparison of the wavelet method with others methods, such as the finite element, finite difference, Adomian decomposition and meshless methods, shows that the proposed wavelet method is accurate and efficient [3-5]. Most interestingly, the corresponding numerical error of the wavelet method is almost independent of the nonlinear intensity of the equations [3, 6]. For the solution of IBVPs using the spatial discretization method, the original equations are usually transformed into a set of ordinary differential equations (ODEs) with unknown time-dependent wavelet coefficients or nodal approximations [7]. The obtained ODEs can be solved either using conventional numerical methods, such as the Runge–Kutta method, or wavelet-based methods. Given the attractive properties of the wavelet method in terms of fast and adaptive computation, studies have been devoted to develop a wavelet-based time discretization scheme for initial value problems (IVPs) [8]. These studies include the vibration of multi-degrees of freedom (MDOF) systems [7-10], Burgers equations [11, 12], nonlinear circuit simulation [13], and nonlinear stiffness systems [14, 15]. However, all these wavelet-based methods treat the initial-valued ODEs as pseudo BVPs [11] by considering a certain time interval. These pseudo BVPs can then be similarly solved by spatial discretization methods, such as the wavelet collocation [7, 10-12] or wavelet Galerkin methods [16]. The disadvantage of such a treatment in terms of the time variable is that the dimension of the resulting algebraic equations depends not only on the dimension of the ODEs after spatial discretization but also on the discretization of the specified time interval [7]. A larger time interval usually results in requiring additional unknown coefficients of wavelet expansion. Therefore, numerical solutions have been investigated only within a significantly limited time range using most of these wavelet methods [9]. One main reason for why using such time discretization strategies has been hypothesized: to control the global error in time domain [7, 11]. However, the simulation of the time response of a complex dynamic system usually requires balance among accuracy, stability, and computation cost [17, 18]. By contrast, the stability of these wavelet time discretization methods are rarely examined, except for the Haar wavelet-based time integration method, which seems limited only to the linear vibration of MDOF systems [9].

---

[1] Corresponding author, Email: jzwang@lzu.edu.cn

In this study, after briefly introducing the properties and construction procedure of the Coiflet, we develop a novel scheme of wavelet approximation of square-integrable functions with an adjustable order of accuracy. Firstly, we use this approximation scheme to propose a high-precision wavelet-based systematic time-integrating scheme for IVPs. Secondly, we also apply the same approximation scheme in the construction of a modified Galerkin method for BVPs. We combine these two methods for the IVPs and BVPs to develop the solution procedure for general nonlinear IBVPs. Numerical examples, including the nonlinear vibration of a Duffing oscillator, Burgers equation [19-23], and Klein Gordon equation [24-27], are considered to demonstrate the accuracy and efficiency of the wavelet-based methods.

**2. Brief introduction to Coiflets**

Coiflets are a class of compactly supported orthogonal wavelets whose wavelet function $\psi(x)$ and scaling function $\varphi(x)$ satisfy the vanishing and shifted vanishing moment properties [28, 29], respectively. More specifically, for the Coiflet with a compact support of [0, 3$N$-1] with $N$ to be a positive even integer, we have [29]:

$$\varphi(x) = \sum_{k=0}^{3N-1} p_k \varphi(2x-k), \tag{1}$$

$$\psi(x) = \sum_{k=0}^{3N-1} (-1)^k p_{3N-1-k} \varphi(2x-k), \tag{2}$$

$$\int_{-\infty}^{\infty} (x-M_1)^n \varphi(x)dx = 0, \text{ and } M_n = M_1^n \text{ for } 1 \leq n \leq N, \tag{3}$$

$$\int_{-\infty}^{\infty} x^n \psi(x)dx = 0, \text{ for } 0 \leq n < N-1, \tag{4}$$

$$\sum_{k\in Z}\varphi(x-k)=1 \tag{5}$$

where $p_k$ is a set of low-pass filter coefficients of the wavelet system, $N$-1 is the number of the vanishing moments, and $M_n = \int_{-\infty}^{\infty} x^n \varphi(x)dx$ is the $n$th-order moment of the scaling function. As shown by Wang [29], these properties can be derived to associate with the filter coefficients $p_k$ as

$$\sum_{k=0}^{3N-1} p_k = 2. \tag{6}$$

$$\sum_{i=0}^{3N-1} p_i p_{i-2k} = 2\delta_{0,k}, \quad k = 0,1,2,\cdots 3N/2 - 1. \tag{7}$$

$$\sum_{j=0}^{3N-1} (-1)^j p_j j^k = 0, \quad k = 0,1,2,\cdots N-1. \tag{8}$$

$$\sum_{j=0}^{3N-1} j^{2i-1} p_j = 2M_1^{2i-1}, \quad i = 1, 2, \ldots, N/2 \tag{9}$$

Eq. (6) provides one equation for the determination of filter coefficients, Eq. (7) provides $3N/2$ equations, Eq. (8) provides $N$-1 independent equations, and Eq. (9) provides $N/2$ equations. Thus, Eqs. (6)–(9) give $3N$ independent equations. Once the $3N$ filter coefficients are obtained, the dyadic values of the scaling, wavelet functions, their derivatives, and integrals can be easily calculated following a standard procedure [29]. Tables 1 and 2 give the values of the scaling function, its derivatives and integrals at integer points, respectively, when $N = 6$ and $M_1 = 7$.

For a function $f(x) \in \mathbf{L}^2(\mathbf{R})$, we denote

$$P^m f(x) = \sum_{k\in\mathbf{Z}} c_{m,k} \varphi_{m,k}(x) \tag{10}$$

where $c_{m,k} = \int_{\mathbf{R}} f(x)\varphi_{m,k}(x)dx$, $\varphi_{m,k}(x) = 2^{m/2}\varphi(2^m x - k)$, and $m$ is sometimes called the resolution level. The sequence of subspaces, $\mathbf{V}_m = \{\sum_k c_{m,k}\varphi_{m,k}(x), m,k \in \mathbf{Z}\}$, satisfy

$$\{0\}\cdots \subset \mathbf{V}_{-2} \subset \mathbf{V}_{-1} \subset \mathbf{V}_0 \subset \mathbf{V}_1 \subset \cdots \subset \mathbf{L}^2(\mathbf{R}), \tag{11}$$

$$P^m f(x) \in \mathbf{V}_m \leftrightarrow P^m f(2x) \in \mathbf{V}_{m+1}, \text{ for all } m\in\mathbf{Z} \text{ and } f(x)\in\mathbf{L}^2(\mathbf{R}), \tag{12}$$

$$P^m f(x) \in \mathbf{V}_m \leftrightarrow P^m f(x+2^{-m}) \in \mathbf{V}_m, \text{ for all } m\in\mathbf{Z} \text{ and } f(x)\in\mathbf{L}^2(\mathbf{R}), \tag{13}$$

$$\lim_{m\to+\infty} \mathbf{V}_m = \bigcup_{m=-\infty}^{+\infty} \mathbf{V}_m \text{ is dense in } \mathbf{L}^2(\mathbf{R}), \tag{14}$$

$$\lim_{m\to-\infty} \mathbf{V}_m = \bigcap_{m=-\infty}^{+\infty} \mathbf{V}_m = \{0\}, \tag{15}$$

The set $\{\phi(x-k), k\in\mathbf{Z}\}$ forms the orthogonal basis of the subspace $\mathbf{V}_0$, (16)

$$x^n \in \mathbf{V}_0, \text{ for } n=1,2,...N\text{-}1. \tag{17}$$

Given that $f(x) \approx P^m f(x)$ and considering Eq. (10), we have

$$f(x) \approx \sum_{k\in\mathbf{Z}} c_{m,k} \varphi_{m,k}(x) \tag{18}$$

The expansion coefficients in Eq. (18) can be calculated as

$$c_{m,k} = \int_{\mathbf{R}} f(x)\varphi_{m,k}(x)dx \approx 2^{-m/2} f(\frac{M_1+k}{2^m}) \qquad (19)$$

which has an accuracy of order $N$ in terms of Eq. (17) [29].

Thus, the function $f(x)$ can be represented as [6]

$$f(x) \approx \tilde{P}^m f(x) = \sum_{k \in Z} f_{M_1+k} \varphi(2^m x - k) \qquad (20)$$

where $f_{M_1+k} = f(x_{M_1+k}), x_{M_1+k} = \frac{M_1+k}{2^m}$, and $f(x) = \tilde{P}^m f(x)$, when $f(x)$ is any polynomial with an order of up to $N$-1, and $\|f(x) - \tilde{P}^m f(x)\|_2 = O(2^{-mN})$, as long as $f(x) \in \mathbf{L}^2(\mathbf{R}) \cap \mathbf{C}^N(\mathbf{R})$.

**Table 1**. Values of the scaling function and its derivatives at integer points for $N = 6$, $M_1 = 7$.

| k | $\varphi(k)$ | $\varphi'(k)$ | $\varphi''(k)$ |
|---|---|---|---|
| 1 | -3.66444437281758E-08 | 7.56061671669701E-07 | 5.85055379740574E-05 |
| 2 | 1.53911669852067E-05 | -1.59556257771372E-04 | -6.23369061345835E-03 |
| 3 | -1.20347428401018E-04 | -4.25387478212341E-04 | 3.77460458441712E-02 |
| 4 | -6.01049457744450E-03 | 3.24199170415698E-02 | 5.03563139649784E-01 |
| 5 | 3.83279869108597E-02 | -2.28004136988034E-01 | -3.47037450030632E+00 |
| 6 | -9.93321926649377E-02 | 9.30506852076968E-01 | 9.76699682300432E+00 |
| 7 | 1.13897129589829E+00 | -2.33731012221323E-01 | -1.33812568105441E+01 |
| 8 | -1.13651812655013E-01 | -5.74622483381206E-01 | 9.06097029996529E+00 |
| 9 | 5.48926083610606E-02 | 7.83291539976493E-02 | -2.93553338773742E+00 |
| 10 | -1.54279915381372E-02 | -4.71311969866472E-03 | 4.40572722586554E-01 |
| 11 | 2.61977141546123E-03 | 5.78218592974164E-04 | -2.37385586325644E-02 |
| 12 | -2.65205890511625E-04 | -1.30140867708115E-04 | 3.01393615774306E-03 |
| 13 | -1.79596295480327E-05 | -4.82963259972682E-05 | 4.13656515185597E-03 |
| 14 | -1.01758785797000E-06 | -7.93241921874722E-07 | 8.38096374634281E-05 |
| 15 | 4.86512575426836E-09 | 2.87068836205473E-08 | -4.90575764840692E-06 |
| 16 | -1.48885946646455E-12 | -1.75812155632504E-11 | 6.01590725201252E-09 |

| k | $\varphi'''(k)$ | $\varphi^{(4)}(k)$ | $\varphi^{(5)}(k)$ |
|---|---|---|---|
| 1 | -7.67106557904889E-06 | 8.74874373281834E-04 | 7.41506226932848E-01 |
| 2 | 4.16578363253693E-04 | -2.46569013532615E-02 | -1.12134042872026E+01 |
| 3 | -9.81913497918669E-03 | 1.56232654203163E-01 | 3.81656638559639E+01 |
| 4 | 3.10644144054426E-03 | 5.35130774147768E-02 | -4.99231511073754E+01 |
| 5 | 4.53761423326168E-01 | -1.24958255797240E+00 | 1.97735417058377E+01 |
| 6 | -5.37009748859017E-01 | 2.50940769618236E+00 | 1.12469149670718E+01 |
| 7 | -1.09291614025045E+00 | -2.68372231247926E+00 | -1.25329879705800E+01 |
| 8 | 2.07047815990482E+00 | 1.96557892020045E+00 | -1.88083891238367E+00 |
| 9 | -8.89472705306856E-01 | -7.77270974961094E-01 | 2.89562961676712E+01 |
| 10 | -5.37514984551947E-02 | -9.30171370312431E-02 | -4.69904951296699E+01 |
| 11 | 6.67430725039445E-02 | 1.75872976647925E-01 | 3.21638529849212E+01 |
| 12 | -1.46891011281847E-02 | -4.24751553744586E-02 | -9.76606877157541E+00 |
| 13 | 3.32547228775262E-03 | 1.00790783057463E-02 | 1.60628120151076E+00 |
| 14 | -1.56035403773163E-04 | -7.74619163766090E-04 | -3.22835615776562E-01 |

| | | |
|---|---|---|
| 15 | -9.13484394985041E-06 | -5.99151722327302E-05 | -2.45202283593051E-02 |
| 16 | 2.24632696100355E-08 | 2.96173316793466E-07 | 2.44913167241909E-04 |

**Table 2.** Values of the scaling function and its integral at integer points for $N = 6$, $M_1 = 7$.

| k | $\varphi^{\int}(k)$ | k | $\varphi^{\int}(k)$ |
|---|---|---|---|
| 1 | -4.527262019294291e-09 | 9 | 9.875787622086671e-01 |
| 2 | 3.793659018256384e-06 | 10 | 1.003591003107596e+00 |
| 3 | 2.929388903777156e-06 | 11 | 9.991452076341547e-01 |
| 4 | -3.134173582887661e-03 | 12 | 1.000073554299164e+00 |
| 5 | 2.264311126637817e-02 | 13 | 1.000001685301915e+00 |
| 6 | -8.486417521139951e-02 | 14 | 1.000000156123390e+00 |
| 7 | 5.234649707003830e-01 | 15 | 9.999999997921172e-01 |
| 8 | 1.051493179855499e+00 | 16 | 1.000000000000944e+00 |

## 3. Approximation of interval-bounded functions

If the function is bounded in the interval $[a, b]$ with $2^m a$ and $2^m b$ as integers, Eq. (20) can be rewritten as

$$f(x) \approx \tilde{P}^m f(x) = \sum_{k=2^m a - 3N + 2 + M_1}^{2^m b + M_1 - 1} f_k \varphi(2^m x - k + M_1) \tag{21}$$

where we still have [6] $\|f(x) - \tilde{P}^m f(x)\|_2 = O(2^{-mN})$. We note that $x_k = k/2^m$ may be located outside the interval $[a, b]$, and the function $f(x)$ has no definitions in these locations. Therefore, the techniques of boundary extension are usually adopted to resolve this problem. In this study, we use power series expansion at each boundary as

$$f(x) = \begin{cases} \sum_{i=0}^{N-1} \dfrac{F_a^{(i)}}{i!}(x-a)^i & x < a \\ f(x) & a \leq x \leq b \\ \sum_{i=0}^{N-1} \dfrac{F_b^{(i)}}{i!}(x-b)^i & x > b \end{cases} \tag{22}$$

where

$$F_b^{(i)} = \begin{cases} f_b^{(i)} = \dfrac{d^i f(b)}{dx^i} & \text{if the } i\text{th derivative of } f(x) \\ & \text{at } x = b \text{ exists and is known} \\ 2^{im} \sum_{k=0}^{\alpha_2} \zeta_{b,i,k} f_{2^m b - k} & \text{otherwise} \end{cases} \tag{23}$$

and similarly,

$$F_a^{(i)} = \begin{cases} f_a^{(i)} = \dfrac{d^i f(a)}{dx^i} & \text{if the } i\text{th derivative of } f(x) \\ & \text{at } x = a \text{ exists and is known} \\ 2^{im} \sum_{k=0}^{\alpha_1} \zeta_{a,i,k} f_{2^m a + k} & \text{otherwise} \end{cases} \tag{24}$$

To derive the coefficients, $\zeta_{b,i,k}$, substituting Eq. (22) into Eq. (21), taking the derivatives on both sides of Eq. (21) with respect to $x$, and considering $x = b$, we have

$$f_b^{(j)} \approx 2^{jm} \sum_{k=2^m b-3N+2+M_1}^{2^m b} f_k \varphi^{(j)}(2^m b - k + M_1)$$
$$+ 2^{jm} \sum_{k=2^m b+1}^{2^m b-1+M_1} \sum_{i=0}^{N-1} \frac{F_b^{(i)}}{i!} (\frac{k}{2^m} - b)^i \varphi^{(j)}(2^m b - k + M_1) \qquad (25)$$

where $m$ is assumed to be sufficiently large so that $2^m b-3N+2+M_1 > 2^m a$. According to Eq. (17), concluding that Eq. (25) is exact for $f(x)$ to be any polynomial with an order of up to $N$-1 when $j < N$ is not difficult. We let the right side of Eq. (25) be equal to $F_b^{(j)}$, and using $F_b^{(j)}$ to replace $f_b^{(j)}$ in the left side of Eq. (25), we have

$$F_b^{(j)} - \sum_{i=0}^{N-1} 2^{jm} F_b^{(i)} \sum_{k=1}^{M_1-1} \frac{1}{i!} (\frac{k}{2^m})^i \varphi^{(j)}(-k+M_1) = 2^{jm} \sum_{k=3N-2-M_1}^{0} f_{2^m b-k} \varphi^{(j)}(k+M_1). \qquad (26)$$

This replacement ensures consistency between the difference-like approximation of the boundary derivatives as shown in the second line of the right side of Eq. (23) and the wavelet approximation of the boundary derivatives in Eq. (25). In addition, Eq. (26) is again exact if $f(x)$ is any polynomial with an order of up to $N$-1, and $F_b^{(j)} = f_b^{(j)}$, for $j=1,2,\ldots,N$-1. Inserting $F_b^{(i)} = 2^{im} \sum_{k=0}^{\alpha_2} \zeta_{b,i,k} f_{2^m b-k}$ into Eq. (26), and taking $\alpha_2 = 3N-2-M_1$, we have

$$\sum_{k=0}^{\alpha_2} [\varphi^{(j)}(k+M_1) - \zeta_{b,j,k} + \sum_{i=0}^{N-1} \zeta_{b,i,k} \sum_{l=1}^{M_1-1} \frac{1}{i!} l^i \varphi^{(j)}(-l+M_1)] f_{2^m b-k} = 0. \qquad (27)$$

Eq. (27) can be written in matrix form as

$$[\mathbf{A}_1 - (\mathbf{I} - \mathbf{B}_1)\mathbf{P}_1]\mathbf{F} = 0 \qquad (28)$$

where $\mathbf{P}_1 = \{\zeta_{b,i,k}\}$, $\mathbf{A}_1 = \{\varphi^{(i)}(k+M_1)\}$, and $\mathbf{B}_1 = \{\sum_{l=1}^{M_1-1} \frac{1}{j!} l^j \varphi^{(i)}(-l+M_1)\}$, $i, j = 0,1,2,\ldots N-1$, $k = 0,1,2,\ldots \alpha_2$.

If Eq. (28) is satisfied for any $\mathbf{F}$, then from Eq. (28), we obtain

$$\mathbf{P}_1 = (\mathbf{I} - \mathbf{B}_1)^{-1} \mathbf{A}_1$$

This result can serve as proof that $f^{(i)}(b) = 2^{im} \sum_{k=0}^{\alpha_2} \zeta_{b,i,k} f_{2^m b-k}$ is exactly satisfied for $f(x) = 1, x, \ldots, x^{N-1}$ (see Appendix A). Similarly, we obtain $\mathbf{P}_0 = \{\zeta_{a,i,k}\}$, $i = 0, 1, 2, \ldots, N$-1, $k=1, 2, \ldots, \alpha_1$, and $\alpha_1 = M_1$-1 (see Appendix B), and $f^{(i)}(a) = 2^{im} \sum_{k=0}^{\alpha_1} \zeta_{a,i,k} f_{2^m a+k}$. Then, substituting Eqs. (23) and (24) into Eq. (22) and taking $x = k/2^m$, respectively, we have

$$f(x) = \begin{cases} \sum_{j=0}^{\alpha_1} f(a + \frac{j}{2^m}) T_{L,j}(x, \boldsymbol{\beta}_L), & x \in (-\infty, 0) \\ f(x), & x \in [a, b] \\ \sum_{j=0}^{\alpha_2} f(b - \frac{j}{2^m}) T_{R,j}(x, \boldsymbol{\beta}_R), & x \in (b, \infty) \end{cases} \qquad (29)$$

in which [30]:

$$T_{L,j}(x,\boldsymbol{\beta}_L) = \sum_{i=0}^{N-1}\beta_{L,i} 2^{im}\frac{\zeta_{a,i,j}}{i!}(x-a)^i, \quad T_{R,j}(x,\boldsymbol{\beta}_R) = \sum_{i=0}^{N-1}\beta_{R,i} 2^{im}\frac{\zeta_{b,i,j}}{i!}(x-b)^i \tag{30}$$

where $\boldsymbol{\beta}_L = \{\beta_{L,i}\}$, $\boldsymbol{\beta}_R = \{\beta_{R,i}\}$, and $i = 0,1,2,3...,N-1$. The introduction of $\boldsymbol{\beta}_L$ and $\boldsymbol{\beta}_R$ merely aims to assign boundary conditions to the function $f(x)$. For example, when $x = a$: $d^i f(a)/dx^i = 0$ and other boundary derivatives are unknown, we simply need to set $\beta_{L,i} = 0$ and all the other elements of $\boldsymbol{\beta}_L$ and $\boldsymbol{\beta}_R$ equal to 1.

Using Eq. (21) to approximate Eq. (29), and making further rearrangement as shown in [30], yields

$$f(x) \approx \sum_{k=2^m a}^{2^m b} f_k \Phi_{m,k}(x) \tag{31}$$

where

$$\Phi_{m,k}(x) = \begin{cases} \varphi(2^m x + M_1 - k) + \sum_{j=2^m a - \alpha_2}^{2^m a} T_{L,k-2^m a}(\frac{j}{2^m},\boldsymbol{\beta}_L)\varphi(2^m x + M_1 - j), & 2^m a \leq k \leq 2^m a + \alpha_1 \\ \varphi(2^m x + M_1 - k), & \alpha_1 + 1 \leq k \leq 2^m b - \alpha_2 - 1, \\ \varphi(2^m x + M_1 - k) + \sum_{j=1+2^m b}^{2^m b + \alpha_1} T_{R,2^m b - k}(\frac{j}{2^m},\boldsymbol{\beta}_R)\varphi(2^m x + M_1 - j), & 2^m b - \alpha_2 \leq k \leq 2^m b \end{cases} \tag{32}$$

and still one has $\left\| f(x) - \sum_{k=2^m a}^{2^m b} f_k \Phi_{m,k}(x) \right\| = O(2^{-mN})$. Given that the boundary extension improves only the smoothness of the boundary regions of the function $g(x)$, the accuracy of approximation is indeed unchanged [6].

As an example, we consider the function $f(x) = e^x$, $x \in [0,1]$ to be approximated by different series expansions such as Fourier series, Chebyshev orthogonal polynomials and the proposed wavelet method by Eq. (31) with $N = 6$, $M_1 = 7$. The absolute errors of these approximations are plotted in Fig. 1 when the number of terms in each series is taken as 16. Fig. 1 shows that approximation using Eq. (31) is much better than the other two methods. The absolute error of the wavelet method of this approximation using the wavelet method is plotted in Fig. 2, which clearly demonstrates that the order of the error is scaled as $2^{-6m}$, that is, the approximation has an order $N$.

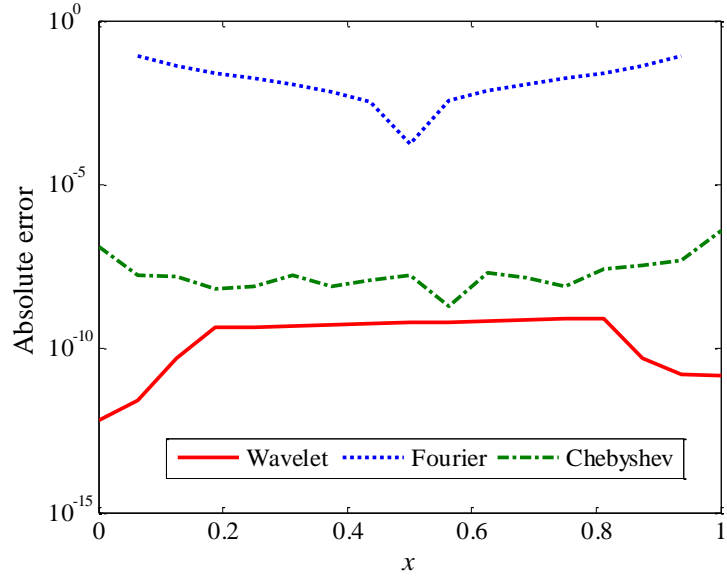

**Fig. 1.** Absolute error of $f(x) = e^x$ using different approximation methods; the scaling function with $N = 6$, $M_1 = 7$ is adopted for the wavelet method.

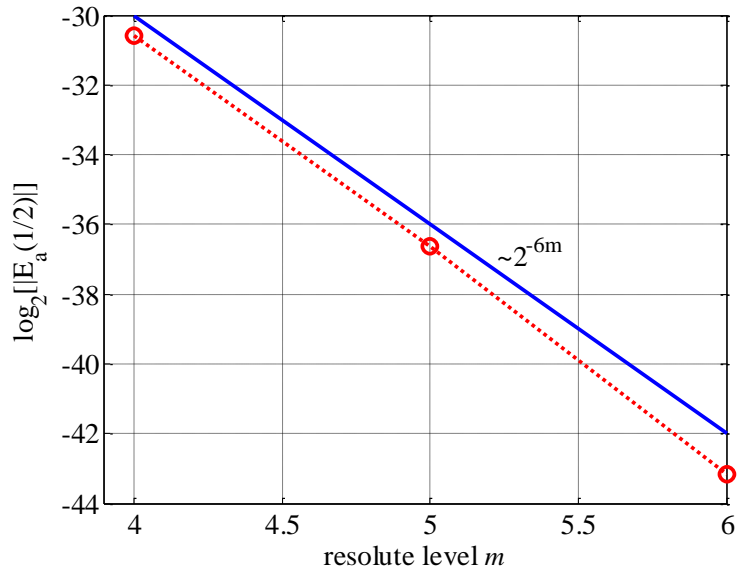

**Fig. 2.** Absolute error, $E_a$, at $x = 1/2$ as a function of the resolute level $m$ in the approximation of $f(x) = e^x$; the scaling function with $N = 6$, $M_1 = 7$ is adopted for the wavelet method.

## 4. Wavelet time-integrating method (WTIM) for the solution of IVPs

In solving IBVPs, one usually needs to spatially discretize the original equations into a set of initial-valued ODEs. By contrast, the ODEs also directly arise from various fields of mathematics and physics, such as nonlinear circuit simulation and vibration of MDOF systems These equations can be generally written as [31]

$$\dot{\mathbf{y}} = \mathbf{f}(t, \mathbf{y}), \quad t \geq 0, \qquad \mathbf{y}_0 = \mathbf{y}(0) \qquad (36)$$

where $\mathbf{y} = [y_1, y_2, \cdots, y_s]^T$ is an unknown vector function of the state variables and $\mathbf{f} = [f_1, f_2, \cdots, f_s]^T$ is a given nonlinear vector function.

To derive the solution method of Eq. (36), we first integrate Eq. (36) as

$$\mathbf{y}(t) = \int_{t-1/2^m}^{t} \mathbf{f}(\tau, \mathbf{y}) d\tau + \mathbf{y}(t - 1/2^m) \tag{37}$$

By setting $a = -\infty, b = j/2^m$ based on Eqs. (31) and (32), we directly approximate every entry of the function vector $\mathbf{f}$ on the time interval $t \in [-\infty, j/2^m]$, $j=1, 2, \ldots$, as

$$f_i(t, \mathbf{y}) \approx \sum_{k=-\infty}^{j} f_i(t_k, \mathbf{y}_k) \Phi_{m,k}(t)$$

where $j$ is a positive integer, $t_k = k/2^m$, and $\mathbf{y}_k = \mathbf{y}(t_k)$. Based on the equation above, we further consider the approximation of $f_i(t, \mathbf{y})$ for $t \in [(j-1)/2^m, j/2^m]$, which can be given by

$$f_i(t, \mathbf{y}) \approx \sum_{k=j-\alpha_2-1}^{j} f_i(t_k, \mathbf{y}_k) \Phi_{m,k}(t) \tag{38}$$

Substituting Eq. (38) into Eq. (37) and taking $t = j/2^m$, we have

$$\mathbf{y}_j \approx \mathbf{y}_{j-1} + \sum_{k=j-\alpha_2-1}^{j} \mathbf{f}(t_k, \mathbf{y}_k) \int_{\frac{j-1}{2^m}}^{\frac{j}{2^m}} \Phi_{m,k}(\tau) d\tau. \tag{39}$$

Rearranging Eq. (39) yields

$$\mathbf{y}_j - \mathbf{f}(t_j, \mathbf{y}_j) \int_{\frac{j-1}{2^m}}^{\frac{j}{2^m}} \Phi_{m,j}(\tau) d\tau \approx \mathbf{y}_{j-1} + \sum_{k=j-\alpha_2-1}^{j-1} \mathbf{f}(t_k, \mathbf{y}_k) \int_{\frac{j-1}{2^m}}^{\frac{j}{2^m}} \Phi_{m,k}(\tau) d\tau. \tag{40}$$

Eq. (40) can be rewritten as

$$\mathbf{y}_j - h\Gamma_{m,0} \mathbf{f}(t_j, \mathbf{y}_j) \approx \mathbf{y}_{j-1} + \sum_{l=1}^{\alpha_2+1} h\Gamma_{m,l} \mathbf{f}(t_{j-l}, \mathbf{y}_{j-l}) \tag{41}$$

where $\Gamma_{m,j-k} = 2^m \int_{\frac{j-1}{2^m}}^{\frac{j}{2^m}} \Phi_{m,k}(\tau) d\tau$ (see Appendix C for the reason why this integral depends only on $m$ and $j-k$.); $\Phi_{m,k}(\tau)$ is given in Eq. (36) with $a = -\infty, b = j/2^m$; and all the elements of $\boldsymbol{\beta}_L$ and $\boldsymbol{\beta}_R$ are equal to 1. In solving Eq. (36) using Eq. (41) at the first $\alpha_2$ steps, we can either set all $\mathbf{y}^{(i)}(0)$, $I = 1, 2, \ldots, N-1$ as unknown or obtain them from Eqs. (36) or (24) in advance. Then, values of $\mathbf{y}_k$, $k = -\alpha_2, \ldots, -1$ in Eq. (41) can be obtained through the power series expansion, $\mathbf{y}(t) \approx \sum_{i=0}^{N-1} \mathbf{y}_0^{(i)} t^i / i!$, $t < 0$.

For the linear ODEs, Eq. (36) can be simplified to

$$\dot{\mathbf{y}} = \mathbf{A}(t)\mathbf{y} + \mathbf{q}(t). \tag{42}$$

Applying Eq. (41) to Eq. (42) yields

$$[\mathbf{I} - h\Gamma_{m,0} \mathbf{A}(t_j)] y_j \approx y_{j-1} + \sum_{l=1}^{\alpha_2+1} h\Gamma_{m,l} [\mathbf{A}(t_{j-l}) y_{j-l} + \mathbf{q}_{j-l}(t_{j-l})] \tag{43}$$

where $\mathbf{I}$ is a unit diagonal matrix. Eqs. (41) and (43) show that the proposed solution method has a an order of accuracy, $N$, which can be any positive even integer, as long as the Coiflet with $N$-1 vanishing moments is adopted.

## 5. Stability analysis of the WTIM

Examining the stability of the WTIM for the solution of IVPs is necessary. We consider the benchmark model problem [32], $\dot{y} = \lambda y$. By applying Eq. (43), we have

$$(1 - \lambda h \Gamma_{m,0}) y_j \approx y_{j-1} + \sum_{k=1}^{\alpha_2+1} \lambda h \Gamma_{m,k} y_{j-k}. \tag{44}$$

Eq. (44) can be rewritten as

$$0 \approx y_{j-1} - y_j + \sum_{k=0}^{\alpha_2+1} \lambda h \Gamma_{m,k} y_{j-k}$$

Furthermore, we have

$$0 \approx -y_j / y_{j-\alpha_2-1} + y_{j-1} / y_{j-\alpha_2-1} + \sum_{k=0}^{\alpha_2+1} \lambda h \Gamma_{m,k} y_{j-k} / y_{j-\alpha_2-1} \tag{45}$$

By defining $\mu = y_{n+1} / y_n$, $z = \lambda h$, according to Eq. (45), we consider

$$z = \frac{p(\mu)}{q(\mu)} \tag{47}$$

where $p(\mu) = \mu^{\alpha_2+1} - \mu^{\alpha_2}$ and $q(\mu) = \sum_{k=0}^{\alpha_2+1} \Gamma_{m,k} \mu^{\alpha_2+1-k}$. Taking Eq. (47) as an algebraic equation of $\mu$ in terms of reference [32], the proposed wavelet method is stable if all the absolute values of the roots of Eq. (47), $|\mu|$, as a function of $z$, are less than 1. The stability domain in terms of $z$ for the method can be determined by investigating the roots of Eq. (46). Letting $\mu = e^{i\theta}$, $0 \leq \theta \leq 2\pi$, we have $|\mu| = 1$. Inserting $\mu = e^{i\theta}$ into Eq. (47), we obtain a closed curve in the $z$-complex plane. The stability domain of the WTIM can be numerically determined, which turns out to be the interior region of the closed curve [32]. As an example, considering the Coiflet with $N = 6$ and $M_1 = 7$, Fig. 3 plots the stability domain enclosed by the blue curve.

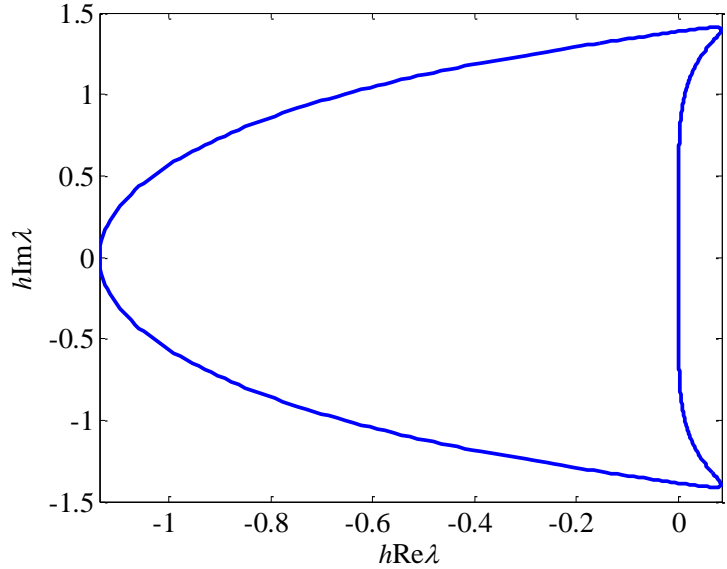

**Fig. 3.** Stability domain (interior of the curves) of the proposed wavelet time-integrating method, where the scaling function with $N = 6$, $M_1 = 7$ is adopted.

## 6. Simultaneous space-time wavelet method (SSTWM) for IBVPs

We consider a class of nonlinear initial-boundary value problems as follows:

$$\frac{\partial u(x,t)}{\partial t} = \mathbf{L}^0 u(x,t) + \mathbf{L}^1 \mathbf{N}[u(x,t), x, t] + f(x,t), \quad 0 \leq x \leq 1, \quad t > 0,$$

$$\partial^{i_s} u / \partial x^{i_s} |_{x=0} = 0, \quad \partial^{j_s} u / \partial x^{j_s} |_{x=1} = 0 \tag{49}$$

$$u(x,0) = g(x)$$

where $u(x,t)$ is the unknown function, $\mathbf{L}^0$ and $\mathbf{L}^1$ denote differential operators, $\mathbf{N}$ is a nonlinear function of $u(x,t)$, and $i_s$, $j_s$, and $s$ are non-negative integers. For the equations with inhomogeneous boundary conditions, a simple transformation can reduce them to that with homogeneous boundary conditions.

For the solution of Eq. (49), we use the proposed wavelet expansion of Eqs. (31) and (32) to approximate various functions in this equation and then adopt the Galerkin method to perform spatial discretization. To satisfy the boundary conditions in Eq. (49), we set $\beta_{L,i_s} = \beta_{L,j_s} = 0$ and all the other elements of $\boldsymbol{\beta}_L$ and $\boldsymbol{\beta}_R$ equal to 1. The resulting modified scaling function basis of $\Phi_{n,k}(x)$ in Eq. (32) is specified accordingly, which is re-denoted as $\tilde{\Phi}_{n,k}(x)$. The functions and nonlinear terms with $x \in [0,1]$ in the Eq. (49) can then be approximated by

$$u(x,t) \approx \sum_{k=0}^{2^n} u_i(x_k,t)\tilde{\Phi}_{n,k}(x), \tag{50}$$

$$\mathbf{N}[u(x,t),x,t] \approx \sum_{k=0}^{2^n} \mathbf{N}[u(x_k,t),x_k,t]\Phi_{n,k}(x), \tag{51}$$

$$f(x,t) \approx \sum_{k=0}^{2^n} f(x_k,t)\Phi_{n,k}(x) \tag{52}$$

where $\Phi_{n,k}(x)$ implies that the boundary derivatives of the corresponding functions are not specified, or the means of all the elements of $\boldsymbol{\beta}_L$ and $\boldsymbol{\beta}_R$ equal to 1.

Substituting Eqs. (50)–(52) into Eq. (49) yields

$$\sum_{k=0}^{2^n} \frac{du(x_k,t)}{dt}\tilde{\Phi}_{n,k}(x) \approx \sum_{k=0}^{2^n} u(x_k,t)\mathbf{L}^0\tilde{\Phi}_{n,k}(x) + \sum_{k=0}^{2^n} \mathbf{N}[u(x_k,t),x_k,t]\mathbf{L}^1\Phi_{n,k}(x) + f(x_k,t)\Phi_{n,k}(x) \tag{53}$$

Applying the Galerkin method to Eq. (53), we then obtain

$$\dot{\mathbf{U}}(t) = \mathbf{A}^{-1}\mathbf{B}\mathbf{U}(t) + \mathbf{A}^{-1}\mathbf{C}\mathbf{V}(t) + \mathbf{A}^{-1}\mathbf{E}\mathbf{F}(t), \quad \mathbf{U}(0) = \mathbf{G} \tag{54}$$

where matrices $\mathbf{A} = \{a_{lk} = \int_0^1 \tilde{\Phi}_{n,k}(x)\tilde{\Phi}_{n,l}(x)dx\}$, $\mathbf{B} = \{b_{lk} = \int_0^1 \mathbf{L}^0\tilde{\Phi}_{n,k}(x)\tilde{\Phi}_{n,l}(x)dx\}$, $\mathbf{C} = \{c_{lk} = \int_0^1 \mathbf{L}^1\Phi_{n,k}(x)\tilde{\Phi}_{n,l}(x)dx\}$, $\mathbf{E} = \{e_{lk} = \int_0^1 \Phi_{n,k}(x)\tilde{\Phi}_{n,l}(x)dx\}$, vectors $\mathbf{U}(t) = \{u_k = u(x_k,t)\}^T$, $\mathbf{V}(t) = \{v_k = \mathbf{N}[u(x_k,t),x_k,t]\}^T$, and $\mathbf{F}(t) = \{f_k = f(x_k,t)\}^T$, $\mathbf{G}(t) = \{g_k = g_k(x)\}$. The entries in these matrices, $a_{lk}$, $b_{lk}$, $c_{lk}$, and $e_{lk}$, can be obtained exactly through the procedure suggested by Wang [29]. Then, Eq. (54) as an IVP can be readily solved using the proposed WTIM. However, most interestingly, we note that in the semi-discretization system (54) of the nonlinear IBVP problem (49), $\mathbf{A}$, $\mathbf{B}$, $\mathbf{C}$, $\mathbf{E}$ and $\mathbf{G}$ are all constant matrices, which are completely independent of the unknown vector $\mathbf{U}(t)$ and time $t$. Thus, in subsequent time integration for solving ordinary differential equations (54), no matrix generated in spatial discretization needs updating, which implies that a fully decoupling between spatial and temporal discretizations is achieved in the present wavelet formulation.

## 7. Numerical examples

Four different dynamic problems are considered to demonstrate the accuracy and efficiency of the proposed wavelet based methods: the free vibration equation of a single oscillator, nonlinear vibration equation of a Duffing oscillator, the Burgers equation, and the Klein–Gordon equation.

**Problem 1:** Free vibration equation of a single oscillator

We consider the free vibration of a single oscillator with the dimensionless governing equation as
$$\ddot{x} + \xi x = 0. \tag{55}$$

We take $\xi = 4\pi^2$, and the initial conditions are $x(0) = 1, \dot{x}(0) = 0$. The analytical solution of Eq. (55) can be immediately obtained as $x(t) = \cos 2\pi t$. We may define $\mathbf{y}(t) = [y_1(t), y_2(t)]^T = [x(t), \dot{x}(t)]^T$; Eq. (55) can then be changed into

$$\dot{\mathbf{y}} = \begin{Bmatrix} 0 & 1 \\ -\xi & 0 \end{Bmatrix} \mathbf{y}, \quad \mathbf{y}(0) = [1, 0]^T \tag{56}$$

Fig. 4 shows the displacement response of this oscillator obtained using the proposed WTIM with Coilfet of $N = 6$ and $M_1 = 7$, the central difference, Newmark-$\beta$, and Wilson-$\theta$ methods, respectively. Fig. 4 shows that the proposed wavelet method has the smallest period elongation among all the classical methods. Fig. 5 plots the absolute error of the WTIM as a function of step size $h$, which shows that WTIM has an accuracy order of 6.5, that is, $N+1/2$. In addition, Fig. 5 shows that the accuracy order of the WTIM is unchanged, whether we set all $\mathbf{y}^{(i)}(0)$, $i = 1,2,...,N-1$ as unknown or obtain them from Eqs. (36) or (24) in advance.

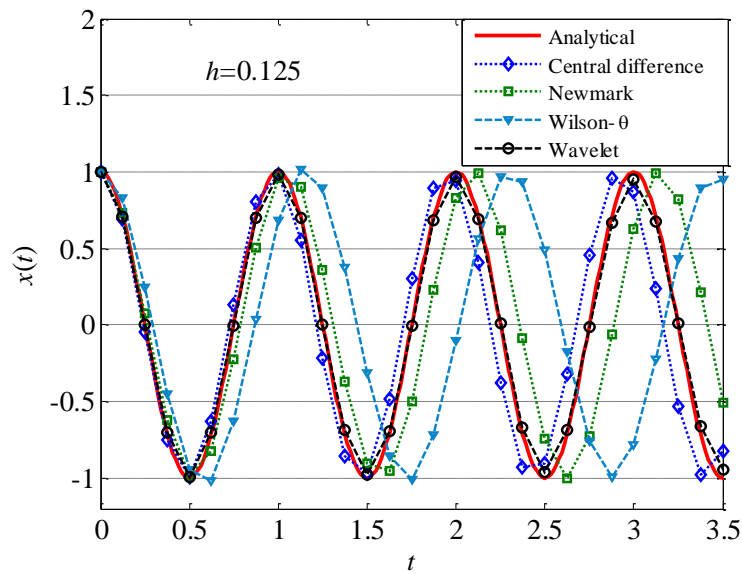

**Fig. 4.** Displacement response of the single oscillator obtained using different numerical methods under the normalized step size $h = 0.125$; the scaling function with $N = 6$, $M_1 = 7$ is adopted for the wavelet method.

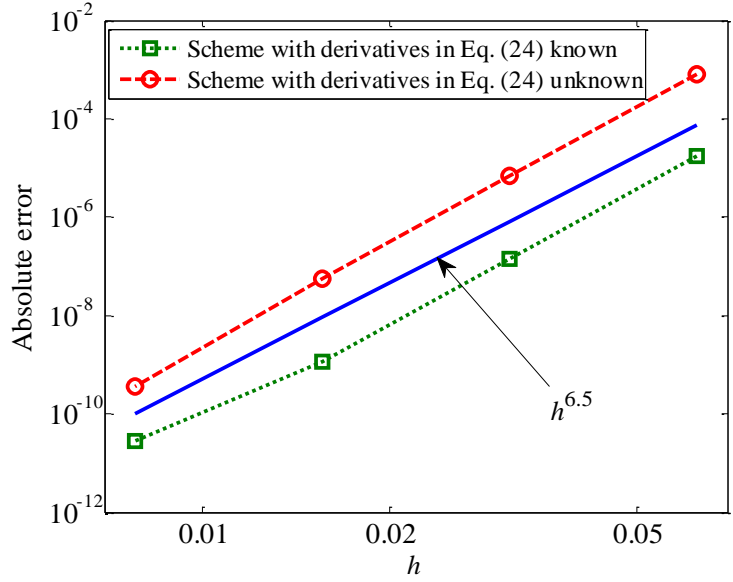

**Fig. 5.** Absolute error on the displacement response at $t = 4$ as a function of time step size $h$, where the scaling function with $N = 6$, $M_1 = 7$ is adopted.

**Problem 2:** Free vibration equation of the Duffing oscillator

The free vibration of the Duffing oscillator represents the vibration of an unforced pendulum with nonlinear restoring force. We choose this problem to test the proposed WTIM because the analytical solution of this nonlinear equation can be obtained by the Jacobi-elliptic functions [18]:

$$\ddot{x} + \omega_n^2 x + \eta x^3 = 0, \quad t \geq 0, \tag{57}$$
$$x(0) = 1, \quad \dot{x}(0) = 0.$$

where $\omega_n = 1$ and the analytical solution is $x(t) = \text{cn}(\sqrt{\eta+1}\, t; k)$ with $k^2 = \eta/(2\eta + 2)$. Denoting $\mathbf{y}(t) = [y_1(t), y_2(t)]^T = [x(t), \dot{x}(t)]^T$, Eq. (57) can be rewritten as

$$\dot{\mathbf{y}} = \begin{bmatrix} 0 & 1 \\ -\omega_n^2 & 0 \end{bmatrix} \mathbf{y} + \begin{bmatrix} 0 & 0 \\ -\eta & 0 \end{bmatrix} \begin{Bmatrix} y_1^3 \\ y_2^3 \end{Bmatrix}, \quad \mathbf{y}(0) = [1, 0]^T. \tag{58}$$

We use the proposed WTIM to solve Eq. (58), where we choose Coiflets with $N = 4$, $M_1 = 7$ and $N = 6$, $M_1 = 7$, respectively. Figs. 6 and 7 show the absolute error of $x(t)$ as a function of step size $h$ at different moments and under different parameters $\eta$, respectively. Figs. 6 and 7 show that the errors decay as $h^{-6}$ for $N = 6$, and decay as $h^{-4}$ for $N = 4$. This example again demonstrates that the accuracy order of this proposed wavelet method is $N$. This accuracy remains unchanged when the intensity of nonlinearity of Eq. (58), characterized by parameter $\eta$, increases.

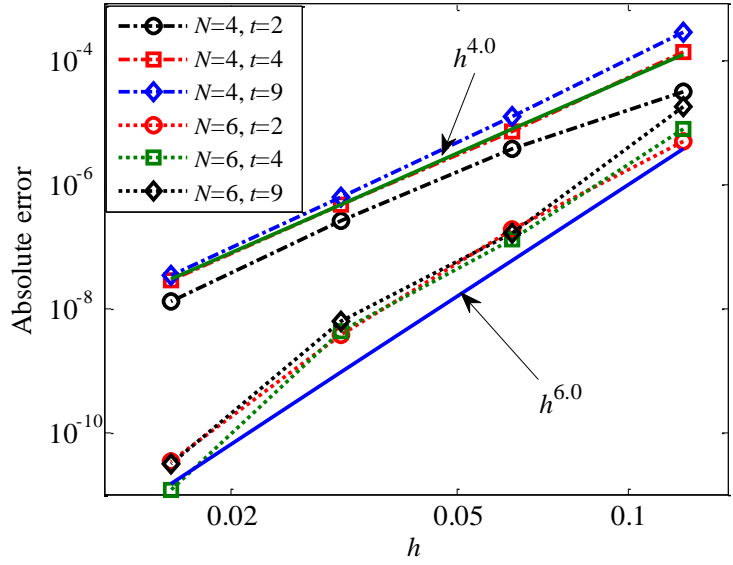

**Fig. 6.** Absolute error of the displacement at different times as a function of time step size $h$ obtained using the Coiflets with $N = 4$, $M_1 = 7$ and $N = 6$, $M_1 = 7$.

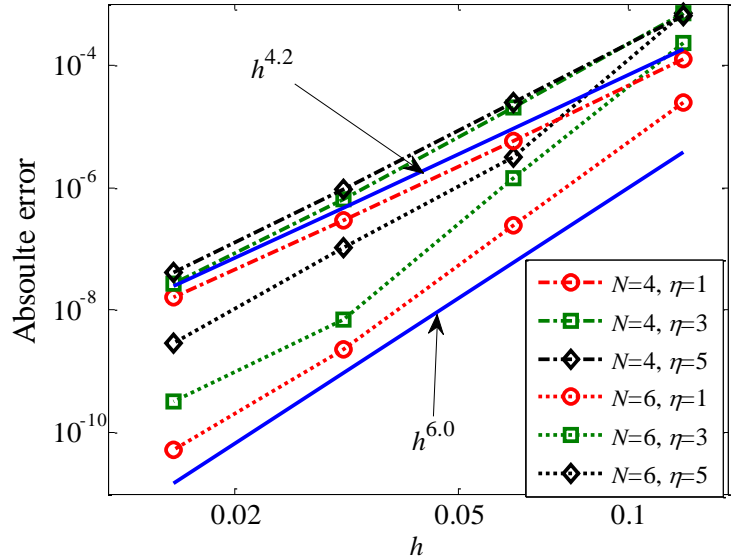

**Fig. 7.** Absolute error of the displacement as a function of time step size $h$ under different parameters $\eta$, obtained using the Coiflets with $N = 4$, $M_1 = 7$ and $N = 6$, $M_1 = 7$.

**Problem 3:** Burgers' equation

Burgers' equation is a fundamental and typical partial differential equation in fluid mechanics and has been widely used to describe wave propagation phenomena in acoustics and hydrodynamics. In this study, we consider the following one-dimensional case:

$$\frac{\partial u}{\partial t} + \frac{1}{2}\frac{\partial u^2}{\partial x} - \frac{1}{\text{Re}}\frac{\partial^2 u}{\partial x^2} = 0, \quad 0 < x < 1, \ t > 0, \tag{61}$$

$$u(0,t) = u(1,t) = 0$$

where the initial boundary condition is taken as $u(x,0) = 2\pi \sin(\pi x)/[100\,\text{Re} + \text{Re}\cos(\pi x)]$ for Case a [19, 23] and $u(x,0) = \sin(\pi x)$ for Case b [20–22]. The exact solutions for these two cases can be found in references [19–23]. Eq. (61) is solved using the proposed simultaneous space-time wavelet method, where the spatial resolution level is taken as $n = 4$. The time step size $h$ and the Reynolds number are taken as 1/64, 200 for Case a and 1/256, 10 for Case b. By contrast, many researchers have solved Eq. (61) using different numerical methods, such as the differential quadrature method [19], methods based on the multiquadratic quasi-interpolation [20], spline interpolation [21], the quadratic B-spline finite element method [22], and the cubic B-splines collocation method [23]. The time derivative is usually discretized using the forward difference methods [19-21] or integration method based on the Runge–Kutta method [23]. Table 3 lists the maximum absolute errors of the numerical solution at $t = 1$ by the proposed wavelet method for $N = 6$, $M_1 = 7$, and several other methods [19–23]. Table 3 shows that that the wavelet method is much more accurate than those proposed in [19–23]. The solutions of Burgers' equation can develop a sharp transition region with significantly small width when the initial condition is smooth [23]. This situation poses a challenge for conventional numerical methods. For the case in which $u(x,0) = \sin(2\pi x)$ and Re = 100, Fig. 8 shows the numerical results on $u(x,t)$ as a function of $x$ at different time $t$, where the spatial resolution level is taken as $n = 6$, and the time step size $h = 1/2^8$. Fig. 8 shows that the numerical solution of Burgers' equation clearly presents a sharp shock wave front after a certain propagation time, implying the competency of the proposed wavelet method in solving such problems.

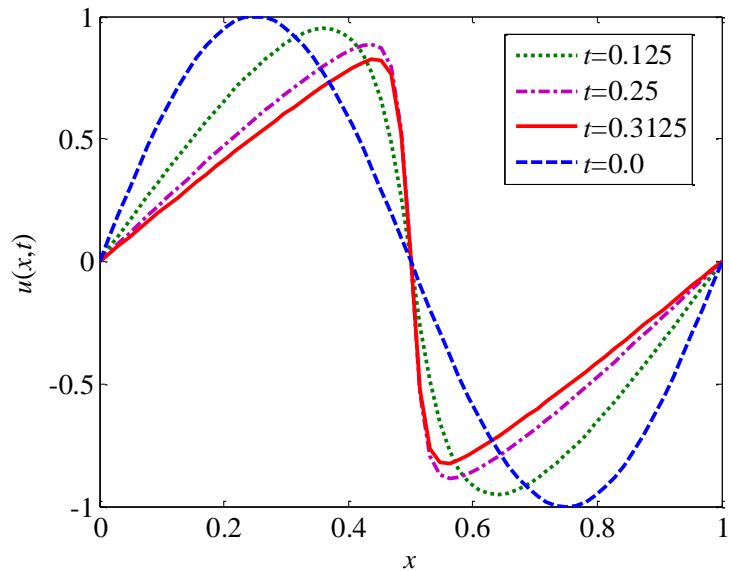

**Fig. 8.** Numerical solutions of Burgers' equation with $u(x,0) = \sin(2\pi x)$, Re = 100. The resolution level of the spatial approximation is $n = 6$, and time step size is $h = 1/2^8$.

**Table 3.** Maximum absolute error of the Burgers equation at $t = 1$ using different methods.

|        | Wavelet    | MCBSCM [19] | WADQM [20] | MQQI [21] | BSQI [22] | FEM [23] |
|--------|------------|-------------|------------|-----------|-----------|----------|
|        | $dx$=1/16  | $dx$=1/20   | $dx$=1/20  | $dx$=1/100| $dx$=1/100| $dx$=1/80|
|        |            | $h$=0.01    | $h$=0.01   | $h$=0.001 | $h$=0.001 | $h$=0.0001|
| case a | 2.0558e-09 | 3.062e-08   | 1.091e-08  | --        | --        | --       |
| case b | 5.8176e-05 | --          | --         | 1.59e-02  | 3.50e-4   | 2.69e-3  |

**Table 4.** Maximum absolute error of the Klein–Gordon equation using different methods, where $dx$ is the spatial "mesh" size, and $h$ the time step size. For the wavelet method, $dx = 1/2^n$, and $n$ is the spatial resolution level.

| $t$ | Wavelet | MFDCM [24] | CM [24] | FMC [25] | RBF [26] | MQ [27] |
|-----|---------|------------|---------|----------|----------|---------|

|     | $dx=0.0625$ $h=1/16$ | $dx=0.0625$ $h=0.0001$ | $dx=0.0625$ $h=0.0001$ | $dx=0.04$ $h=0.001$ | $dx=0.02$ $h=0.0001$ | $dx=0.025$ $h=0.0001$ |
|-----|----------------------|------------------------|------------------------|---------------------|----------------------|-----------------------|
| 0.5 | 4.144e-09            | 5.6e-07                | 4.5e-05                | --                  | --                   | --                    |
| 1.0 | 3.430e-08            | 8.2e-07                | 4.6e-05                | 2.780e-04           | 1.101e-05            | 3.176e-04             |
| 2.0 | 2.716e-07            | --                     | --                     | 2.981e-05           | 1.650e-04            | --                    |
| 3.0 | 8.764e-07            |                        | --                     | 1.469e-04           | 5.973e-04            |                       |
| 4.0 | 1.944e-06            | --                     | --                     | 1.165e-04           | 1.826e-03            | --                    |
| 5.0 | 3.489e-06            |                        | --                     | 4.306e-05           | 3.692e-03            |                       |

**Problem 4:** Klein–Gordon equation

The nonlinear Klein–Gordon equation has an important role in mathematical physics, particularly in investigating the behavior and interaction of solutions in condensed matter physics, nonlinear wave equations, and dispersive phenomena in relativistic physics. A one-dimensional case can be given by

$$\frac{\partial^2 u}{\partial t^2} - \frac{\partial^2 u}{\partial x^2} + u^2 = f(x,t), 0 < x < 1, \ t > 0, \qquad (62)$$

$$u(x,0) = 0, \ u_t(x,0) = 0.$$

When $f(x,t) = 6xt(x^2-t^2)+x^6 t^6$, the analytical solution [24-27] is $x(t) = x^3 t^3$. Eq. (26) has been solved by the mixed finite difference and collocation methods [24-27], where time discretization is implemented using the finite difference [24, 25, 27] or Runge–Kutta method [26], and spatial discretization is implemented using the collocation method associated with the spline basis [24-25] and the radial basis functions [26–27], respectively. To use the proposed wavelet method, we define

$$w_1(t) = u(t), w_2(t) = \frac{\partial u}{\partial t} \qquad (63)$$

Eq. (26) can then be changed into

$$\frac{\partial w_1}{\partial t} = w_2$$
$$\frac{\partial w_2}{\partial t} - \frac{\partial^2 w_1}{\partial x^2} + w_1^2 = f(x,t) \qquad (64)$$

Taking the resolution level as $n = 4$ in the space domain, Eq. (64) can be discretized into a system of first-order nonlinear ODEs with 34 degrees of freedom. The maximum absolute errors of the numerical solution obtained by the proposed wavelet method and the three other methods [24-27] are given in Table 4. The proposed SSTWM has much better accuracy even under a larger time step than that of the finite difference collocation, finite element collocation, and radial basis function meshless methods.

**5. Conclusions**

Based on the unique property of the Coiflets in the approximation of functions, we propose a procedure of a space–time fully decoupled discretization method of high precision to solve nonlinear IBVPs. Using the wavelet Galerkin method, the PDEs are usually transformed into a set of ODEs with unknown nodal parameters that depend on time. We then construct a novel wavelet-based step-by-step time-integrating scheme for ODEs, which has $N$th-order accuracy. Most promisingly, this order of accuracy can be adjusted to any even number when the Coiflet with the order of vanishing moments of this number minus one is adopted. Different numerical examples show that the proposed wavelet-based solution method has much better accuracy and efficiency than most existing numerical methods. We expect that the proposed wavelet method to have promising applications in solving a broad range of nonlinear dynamic problems in condensed matter physics.

**Appendix A**

We first consider the following relation:

$$\mathbf{A}_1 \mathbf{F} = \left\{ \begin{array}{c} \sum_{k=0}^{\alpha_2} \varphi(k+M_1) f_{2^m b-k} \\ \ldots \\ \sum_{k=0}^{\alpha_2} \varphi^{(N-1)}(k+M_1) f_{2^m b-k} \end{array} \right\} = \left\{ \begin{array}{c} \sum_{l=-\alpha_2}^{0} (\frac{2^m b + l}{2^m})^{N-1} \varphi(-l+M_1) \\ \ldots \\ \sum_{l=-\alpha_2}^{0} (\frac{2^m b + l}{2^m})^{N-1} \varphi^{(N-1)}(-l+M_1) \end{array} \right\}. \tag{A1}$$

We then consider

$$(\mathbf{I}-\mathbf{B}_1) \left\{ \begin{array}{c} f(b) \\ \ldots \\ 2^{-im} f^{(i)}(b) \\ \ldots \\ 2^{-(N-1)m} f^{(N-1)}(b) \end{array} \right\} = \left\{ \begin{array}{c} f(b) - \sum_{j=0}^{N-1} \sum_{l=1}^{M_1-1} f^{(j)}(b) \frac{l^j}{2^{mj} j!} \varphi(-l+M_1) \\ \ldots \\ 2^{-im} f^{(i)}(b) - \sum_{j=0}^{N-1} \sum_{l=1}^{M_1-1} f^{(j)}(b) \frac{l^j}{2^{mj} j!} \varphi^{(i)}(-l+M_1) \\ \ldots \\ 2^{-(N-1)m} f^{(N-1)}(b) - \sum_{j=0}^{N-1} \sum_{l=1}^{M_1-1} f^{(j)}(b) \frac{l^j}{2^{mj} j!} \varphi^{(N-1)}(-l+M_1) \end{array} \right\}. \tag{A2}$$

For $f(x) = x^n$, $n = 0, 1, 2, \ldots, N-1$, Eq. (A2) becomes

$$(\mathbf{I}-\mathbf{B}_1) \left\{ \begin{array}{c} f(b) \\ \ldots \\ 2^{-im} f^{(i)}(b) \\ \ldots \\ 2^{-(N-1)m} f^{(N-1)}(b) \end{array} \right\} = \left\{ \begin{array}{c} f(b) - \sum_{j=0}^{n} \sum_{l=1}^{M_1-1} C_j^n b^{n-j} (\frac{l}{2^m})^j \varphi(-l+M_1) \\ \ldots \\ 2^{-im} f^{(i)}(b) - \sum_{j=0}^{n} \sum_{l=1}^{M_1-1} C_j^n b^{n-j} (\frac{l}{2^m})^j \varphi^{(i)}(-l+M_1) \\ \ldots \\ 2^{-(N-1)m} f^{(N-1)}(b) - \sum_{j=0}^{n} \sum_{l=1}^{M_1-1} C_j^n b^{n-j} (\frac{l}{2^m})^j \varphi^{(N-1)}(-l+M_1) \end{array} \right\}$$

$$(\mathbf{I}-\mathbf{B}_1) \left\{ \begin{array}{c} f(b) \\ \ldots \\ 2^{-im} f^{(i)}(b) \\ \ldots \\ 2^{-(N-1)m} f^{(N-1)}(b) \end{array} \right\} = \left\{ \begin{array}{c} f(b) - \sum_{l=1}^{M_1-1} (b+\frac{l}{2^m})^n \varphi(-l+M_1) \\ \ldots \\ 2^{-im} f^{(i)}(b) - \sum_{l=1}^{M_1-1} (b+\frac{l}{2^m})^n \varphi^{(i)}(-l+M_1) \\ \ldots \\ 2^{-(N-1)m} f^{(N-1)}(b) - \sum_{l=1}^{M_1-1} (b+\frac{l}{2^m})^n \varphi^{(N-1)}(-l+M_1) \end{array} \right\} \tag{A3}$$

.

For $f(x) = x^n$, $i, n = 0, 1, 2,\ldots, N-1$, considering Eq. (21) and noting Eq. (17), we have

$$f^{(i)}(x) = 2^{im} \sum_{k=2^m a-3N+2+M_1}^{2^m b+M_1-1} f_k \varphi^{(i)}(2^m x - k + M_1). \tag{A4}$$

Inserting $x = b$ yields

$$f^{(i)}(b) = \sum_{k=2^m b-3N+2+M_1}^{2^m b+M_1-1} f(\frac{k}{2^m})\varphi^{(i)}(2^m b - k + M_1) \tag{A5}$$

or

$$2^{-im} f^{(i)}(b) - \sum_{l=1}^{M_1-1} f(b+\frac{l}{2^m})\varphi^{(i)}(2^m b - k + M_1) = \sum_{l=-\alpha_2}^{0} f(b+\frac{l}{2^m})\varphi^{(i)}(2^m b - k + M_1). \tag{A6}$$

From Eqs. (A1), (A3), and (A6), we can derive

$$\mathbf{A}_1 \mathbf{F} = (\mathbf{I} - \mathbf{B}_1) \mathbf{F}_d \tag{A7}$$

where $\mathbf{F}_d = \{2^{-im} f^{(i)}(b)\}$, $f(x) = x^n$, $i, n = 0, 1, 2,\ldots N-1$. From Eq. (A7), we can further obtain

$$\mathbf{F}_d = (\mathbf{I} - \mathbf{B}_1) \mathbf{A}_1 \mathbf{F} \tag{A8}$$

or

$$f^{(i)}(b) = 2^{im} \sum_{k=0}^{\alpha_2} \zeta_{b,i,k} f_{2^m b-k} \quad , i = 0, 1, 2,\ldots, N-1. \tag{A9}$$

Then, for $f(x) = x^n$, $n = 0, 1, 2,\ldots N-1$, if we define $\mathbf{P}_1 = \{\zeta_{b,i,k}\} = (\mathbf{I} - \mathbf{B}_1)^{-1} \mathbf{A}_1$, Eq. (A9) is exact.

**Appendix B**

To derive the coefficients, $\zeta_{a,i,k}$, substituting Eq. (23) into Eq. (22), taking the derivatives on both sides of Eq. (10) with respect to $x$, and then taking $x = a$, we have

$$\begin{aligned} f_a^{(j)} &\approx 2^{jm} \sum_{k=2^m a}^{2^m a+M_1-1} f_k \varphi^{(j)}(2^m a - k + M_1) \\ &+ 2^{jm} \sum_{k=2^m a-3N+2+M_1}^{2^m a-1} \sum_{i=0}^{N-1} \frac{F_a^{(i)}}{i!}(\frac{k}{2^m} - a)^i \varphi^{(j)}(2^m a - k + M_1) \end{aligned} \tag{B1}$$

where we assume that $m$ is sufficiently large so that $2^m a+M_1-1<2^m b$. Using $F_a^{(j)}$ to replace $f_a^{(j)}$ in Eq. (B1), we have

$$F_a^{(j)} - \sum_{i=0}^{N-1} 2^{jm} F_a^{(i)} \sum_{k=-\alpha_2}^{-1} \frac{1}{i!}(\frac{k}{2^m})^i \varphi^{(j)}(-k + M_1) \approx 2^{jm} \sum_{k=0}^{M_1-1} f_{2^m a+k} \varphi^{(j)}(M_1 - k). \tag{B2}$$

Inserting $F_a^{(i)} = 2^{im} \sum_{k=0}^{\alpha_1} \zeta_{a,i,k} f_{2^m a+k}$ into Eq. (B2), and taking $\alpha_1 = M_1 - 1$, we have

$$\sum_{k=0}^{\alpha_1} [\varphi^{(j)}(M_1 - k) - \zeta_{a,j,k}] f_{2^m a+k} + \sum_{i=0}^{N-1}\sum_{k=0}^{\alpha_1} \zeta_{a,i,k} [\sum_{l=-\alpha_1}^{-1} \frac{1}{i!} l^i \varphi^{(j)}(-l + M_1)] f_{2^m a+k} \approx 0. \tag{B3}$$

Eq. (B3) can be written in the matrix form of $\mathbf{P}_0 = \{\zeta_{a,i,k}\} = (\mathbf{I} - \mathbf{B}_0)^{-1} \mathbf{A}_0$, and matrices

$$\mathbf{A}_0 = \{\varphi^{(i)}(M_1 - k)\} \text{ and } \mathbf{B}_0 = \{\sum_{l=-\alpha_1}^{-1} \frac{1}{i!} l^i \varphi^{(j)}(-l + M_1)\}, i = 0,1,2,\ldots N-1, k = 0,1,2,\ldots \alpha_1. \tag{B4}$$

**Appendix C**

In the case of $a = -\infty$ and $b = j/2^m$, Eq. (32) becomes:

$$\Phi_{m,k}(x) = \begin{cases} \varphi(2^m x + M_1 - k), & k \le j - \alpha_2 \\ \varphi(2^m x + M_1 - k) + \sum_{i=1+j}^{j+\alpha_1} T_{R,j-k}(\frac{i}{2^m}, \boldsymbol{\beta}_R)\varphi(2^m x + M_1 - i), & j - \alpha_2 \le k \le j \end{cases} \quad (C1)$$

We then obtain

$$2^m \int_{(j-1)/2^m}^{j/2^m} \Phi_{m,k}(x)dx$$
$$= \begin{cases} \varphi^{\int}(j-k+M_1) - \varphi^{\int}(j-k+M_1-1), & k \le j-\alpha_2 \\ \varphi^{\int}(j-k+M_1) - \varphi^{\int}(j-k+M_1-1) + \sum_{l=1}^{\alpha_1} T_{R,j-k}(\frac{j+l}{2^m}, \boldsymbol{\beta}_R)[\varphi^{\int}(M_1-l) - \varphi^{\int}(M_1-1-l)], & j-\alpha_2 \le k \le j \end{cases} \quad (C2)$$

We define

$$\Lambda_{m,j-k,l} = T_{R,j-k}(\frac{j+l}{2^m}, \boldsymbol{\beta}_R) = \sum_{i=0}^{N-1} \beta_{R,i} 2^{im} \frac{\zeta_{b,i,j-k}}{i!} (\frac{j+l}{2^m} - \frac{j}{2^m})^i = \sum_{i=0}^{N-1} \beta_{R,i} 2^{im} \frac{\zeta_{b,i,j-k}}{i!} (\frac{l}{2^m})^i \quad (C3)$$

Eq. (C2) can then be rewritten as

$$2^m \int_{(j-1)/2^m}^{j/2^m} \Phi_{m,k}(x)dx = \begin{cases} \varphi^{\int}(j-k+M_1) - \varphi^{\int}(j-k+M_1-1), & \alpha_2 \le j-k \\ \varphi^{\int}(j-k+M_1) - \varphi^{\int}(j-k+M_1-1) + \sum_{l=1}^{\alpha_1} \Lambda_{m,j-k,l}[\varphi^{\int}(M_1-l) - \varphi^{\int}(M_1-1-l)], & 0 \le j-k < \alpha_2 \end{cases} \quad (C4)$$

Eq. (C4) shows that for any $j$ and $k$, the integral, $2^m \int_{(j-1)/2^m}^{j/2^m} \Phi_{m,k}(x)dx$, depends only on $j-k$; therefore, we have

$$\Gamma_{m,j-k} = 2^m \int_{(j-1)/2^m}^{j/2^m} \Phi_{m,k}(x)dx.$$

**Acknowledgement**

This research is supported by grants from the National Natural Science Foundation of China (11472119, 11121202), and the National Key Project of Magneto-Constrained Fusion Energy Development Program (2013GB110002).